\documentclass{article}

\usepackage{amsfonts}
\usepackage{amssymb}
\usepackage{amsthm}
\usepackage{amsmath}
\usepackage{newlfont}
\usepackage{graphicx}
\usepackage[all,arc]{xy}
\usepackage{epsfig}

\usepackage{multicol}
\usepackage[T1]{fontenc}

\newtheorem {theorem}{Theorem}[section]

\newtheorem {corollary}{Corollary}[section]

\author{{\DJ{}or\dj{}e Barali\' {c}}\\ {\small Mathematical Institute SASA}\\[-2mm] {\small Belgrade, Serbia} }

\title{Around the Carnot theorem}
\date{}
\begin{document}
\maketitle

\begin{abstract} We study the Carnot theorem and the configuration of points and lines in connection with it.
It is proven that certain significant points in the configuration
lie on the same lines and same conics. The proof of an equivalent
statement formulated by Bradley is given. An open conjecture,
established by Bradley, is proved using the theorems of Carnot and
Menelaus.
\end{abstract}

\renewcommand{\thefootnote}{}
\footnotetext{This research was supported by the Grant 174020 of
the Ministry for Education and Science of the Republic of Serbia
and Project Math Alive of the Center for Promotion of Science,
Serbia and Mathematical Institute SASA.}

\section{Introduction}

Carnot's theorem can be considered as a generalization of Ceva's
theorem. The theorem of Carnot gives a necessary and sufficient
condition for two points on each side of a triangle to form a
conic.

\begin{theorem}[Carnot's theorem] Let $\triangle A B C$ be a triangle and let $A_1$, $A_2$ be the
points on the line $B C$, $B_1$, $B_2$ on the line $C A$ and $C_1$
and $C_2$ on the line $A B$. The points $A_1$, $A_2$, $B_1$,
$B_2$, $C_1$ and $C_2$ lie on the same conic $\mathcal{C}$ if and
only if
\begin{equation}\label{r}\frac{\overrightarrow{A C_1}}{\overrightarrow{C_1 B}} \cdot
\frac{\overrightarrow{A C_2}}{\overrightarrow{C_2 B}} \cdot
\frac{\overrightarrow{B A_1}}{\overrightarrow{A_1 C}} \cdot
\frac{\overrightarrow{B A_2}}{\overrightarrow{A_2 C}}\cdot
\frac{\overrightarrow{C B_1}}{\overrightarrow{B_1 A}}\cdot
\frac{\overrightarrow{C B_2}}{\overrightarrow{B_2
A}}=1.\end{equation}
\end{theorem}

\medskip

In Section 2 we give a classical proof of Carnot's theorem, using
the theorems of Menelaus and Pascal. This proof can be found in
\cite{Hatt}. We also study some natural points and lines involved
in the configuration and its relations to the side lines of
triangle $\triangle A B C$. Theorems \ref{t2} and \ref{t5}
summarize these results. These theorems are generalizations of
classical Euclidean theorems for incircle of a triangle.

\begin{figure}[h!h!h!]
\centerline{\includegraphics[width=0.8\textwidth]{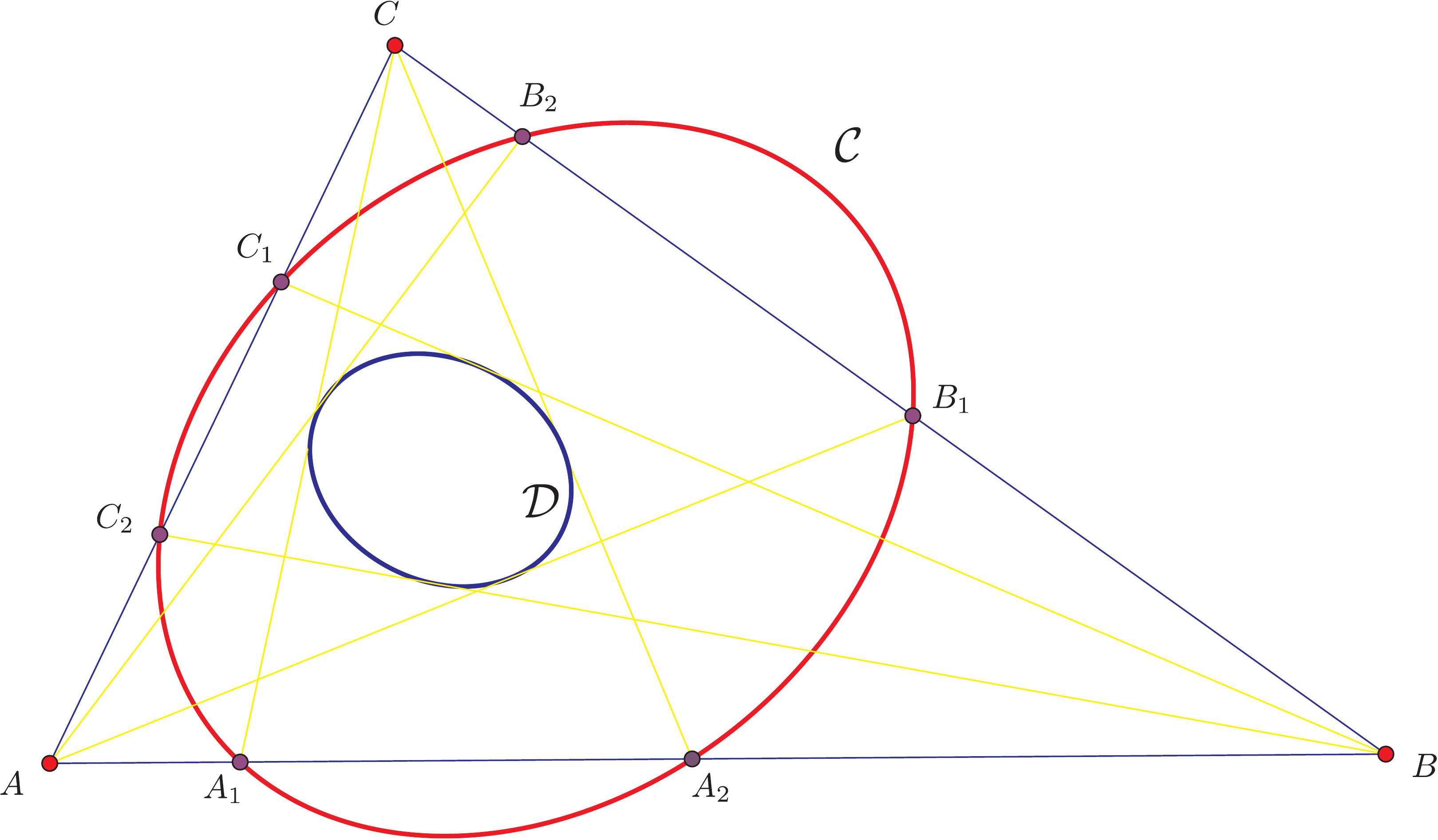}}
\caption{Bradley's theorem} \label{bradley}
\end{figure}

In Section 3 we give an synthetic proof of the following statement
(see Figure \ref{bradley}) which was the first time formulated in
\cite{Bradley}:

\begin{theorem}[Bradley's theorem]\label{glavna} There is a conic $\mathcal{D}$ such that the lines $A
A_1$, $A A_2$, $B B_1$, $B B_2$, $C C_1$ and $C C_2$ are tangents
of $\mathcal{D}$ if and only if the points $A_1$, $A_2$, $B_1$,
$B_2$, $C_1$ and $C_2$ lie on the same conic $\mathcal{C}$.
\end{theorem}

Our goal is to prove an equivalent statement, Corollary
\ref{bradley1} which together with the Poncelet Triangle theorem
implies Bradley's theorem.

\begin{figure}[h!h!h!]
\centerline{\includegraphics[width=0.8\textwidth]{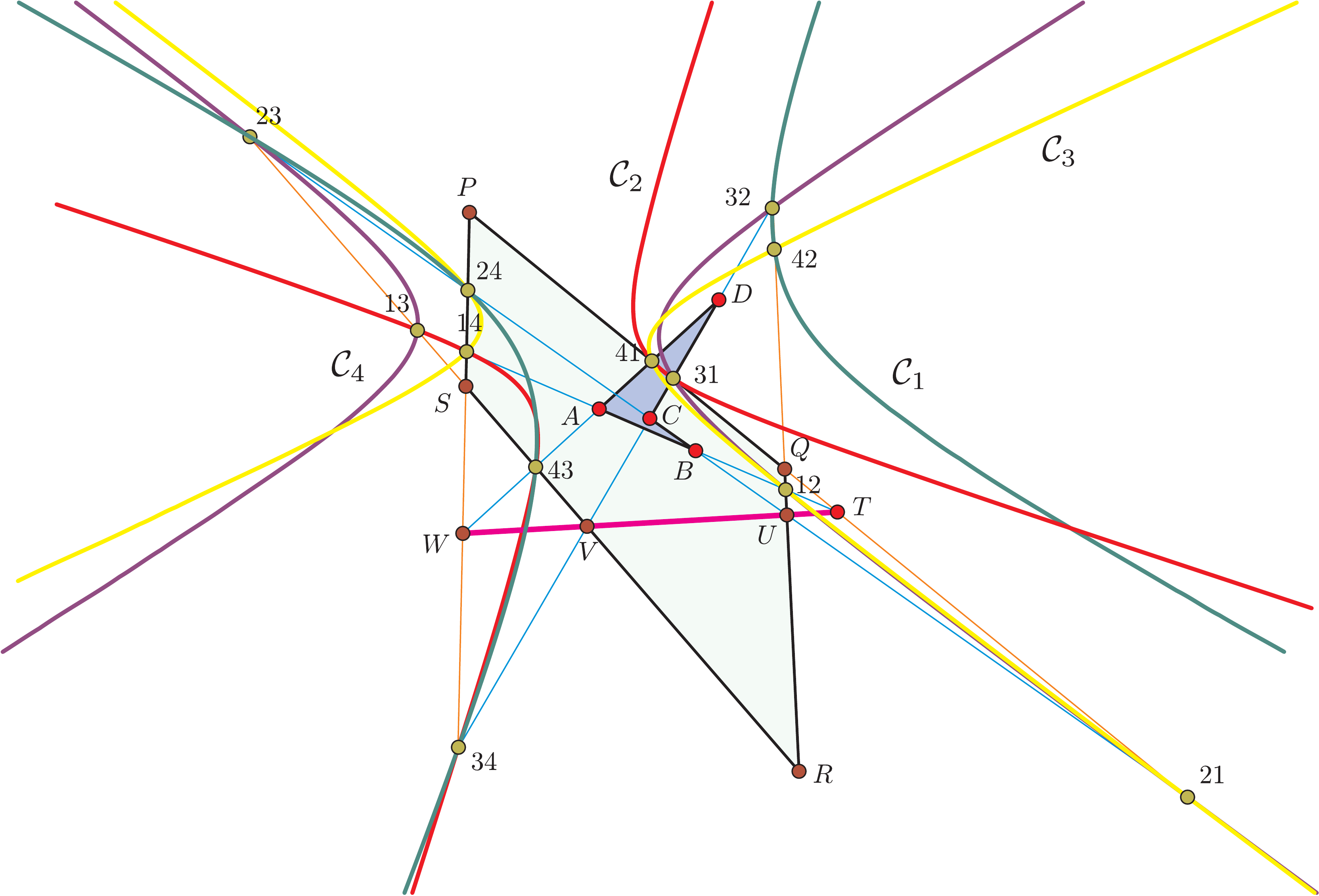}}
\caption{Bradley's theorem about quadrilaterals} \label{bradleygs}
\end{figure}

In the paper \cite{Bradley}, Bradley formulated the following
conjecture (see Figure \ref{bradleygs}):

\begin{theorem}[Bradley's theorem about quadrilaterals]\label{glavnag} Let $ A B C D$ and $P Q R S$ be quadrilateral which are
in axial perspective, that is $T= A B\cap PQ$, $U = BC\cap QR$, $V
= CD\cap S$, $W = DA\cap SP$ are collinear. The other twelve
intersections of the sides of the quadrilaterals are marked with
notation exemplified by $13 = AB\cap RS$, $42 = DA\cap QR$ etc, in
such way that number $1$ corresponds to the sides $A B$ and $P Q$,
$2$ to $B C$ and $Q R$, $3$ to $C D$ and $R S$ and $4$ to $D A$
and $S P$. Then there exist four conics $\mathcal{C}_1$,
$\mathcal{C}_2$, $\mathcal{C}_3$ and $\mathcal{C}_4$ such that the
points $23$, $24$, $32$, $34$, $42$, $43$ lie on conic
$\mathcal{C}_1$,  the points $13$, $14$, $31$, $34$, $42$, $43$
lie on conic $\mathcal{C}_2$, the points $12$, $14$, $21$, $24$,
$41$, $42$ lie on conic $\mathcal{C}_3$ and $12$, $13$, $21$,
$23$, $31$, $32$ lie on conic $\mathcal{C}_4$.
\end{theorem}

This theorem is proved in Section 4.

\medskip

Theorems of Ceva, Menelaos and Carnot are used in \cite{JRG} as
'prototheorems' to build new theorems that involve lines and
conics. It is shown in \cite{JRG} and \cite{Gerb} that any
oriented triangulated 2-manifold can be a frame. This procedure
works for theorems studied in this paper as well. Deep relation
among classical projective geometry and more advanced topics in
mathematics and computer science is explained in
\textit{Perspectives on Projective Geometry}, an inspirative book
by J\"{u}rgen Richter-Gerbert, \cite{Gerb}. Software 'Cinderella'
developed by Ulrich Kortenkamp and  J\"{u}rgen Richter-Gerbert is
used as experimental tool for discovering new results about
Carnot's configuration.

\section{Carnot's theorem}

We start this section with proof of the Carnot theorem.

\begin{figure}[!h!h]
\centerline{\includegraphics[width=0.5\textwidth]{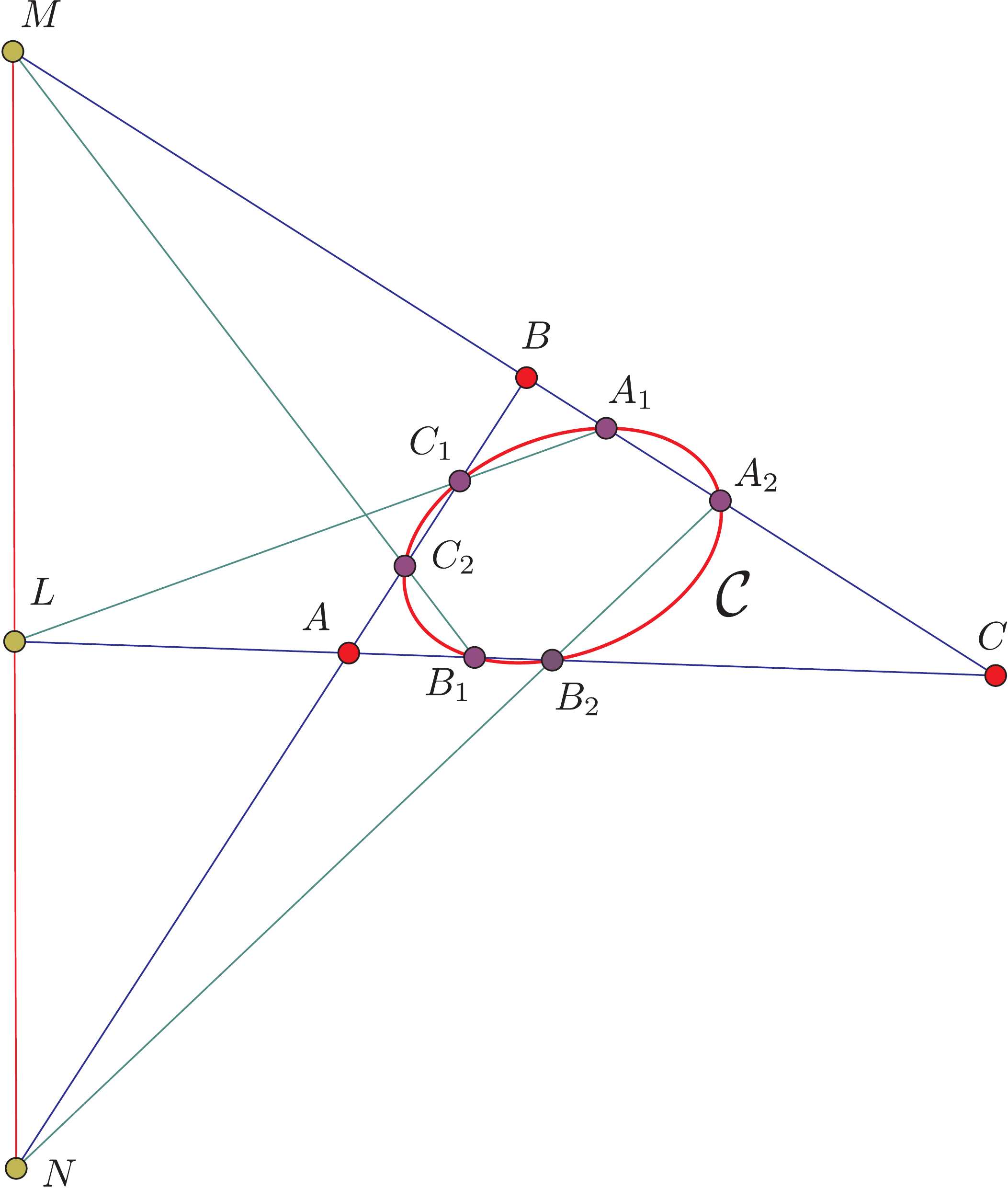}}
\caption{The Carnot theorem} \label{carnot1}
\end{figure}

\noindent \textbf{{Proof of Carnot's theorem:}} Let the points
$A_1$, $A_2$, $B_1$, $B_2$, $C_1$ and $C_2$ lie on the same conic
$\mathcal{C}$ and let $L$ be the intersection of the lines $A_1
C_1$ and $A C$, $M$ the intersection of the lines $B_1 C_2$ and $B
C$ and $N$ the intersection of the lines $A_2 B_2$ and $A B$,
Figure \ref{carnot1}. By the Pascal theorem, the points $L$, $M$
and $N$ lie on the same line, and from the Menelaos theorem the
following holds:\begin{equation}\label{r1} \frac{\overrightarrow{A
L}}{\overrightarrow{L C}}\cdot \frac{\overrightarrow{C
M}}{\overrightarrow{M B}}\cdot \frac{\overrightarrow{B
N}}{\overrightarrow{N A}}=-1. \end{equation}

Applying the Menelaos theorem three times for the lines $A_1 C_1$,
$B_1 C_2$ and $A_2 B_2$ and $\triangle A B C$, we obtain:
\begin{equation}\label{r2}
 \fbox{\parbox{7 mm}{\Large{$\frac{\overrightarrow{A
L}}{\overrightarrow{L C}}$}}}\cdot \frac{\overrightarrow{C
A_1}}{\overrightarrow{A_1 B}}\cdot \frac{\overrightarrow{B
C_1}}{\overrightarrow{C_1 A}}=-1,\end{equation}
\begin{equation}\label{r3} \frac{\overrightarrow{A
B_1}}{\overrightarrow{B_1 C}}\cdot \fbox{\parbox{7
mm}{\Large{$\frac{\overrightarrow{C M}}{\overrightarrow{M
B}}$}}}\cdot \frac{\overrightarrow{B C_2}}{\overrightarrow{C_2
A}}=-1,
\end{equation}
\begin{equation}\label{r4} \frac{\overrightarrow{A
B_2}}{\overrightarrow{B_2 C}}\cdot \frac{\overrightarrow{C
A_2}}{\overrightarrow{A_2 B}}\cdot \fbox{\parbox{7
mm}{\Large{$\frac{\overrightarrow{B N}}{\overrightarrow{N
A}}$}}}=-1.
\end{equation}

Multiplying the relations (\ref{r2}), (\ref{r3}) and (\ref{r4})
and division by (\ref{r1}), yields the relation (\ref{r}).

In the opposite direction, the proof is similar. By the Menelaos
theorem, the relations (\ref{r2}), (\ref{r3}) and (\ref{r4}) hold.
From the relations (\ref{r2}), (\ref{r3}), (\ref{r4}) and
(\ref{r}) one can easily deduce the relation (\ref{r1}), so by the
converse of the Menelaos theorem, the points $L$, $M$ and $N$ lie
on the same line. The converse of the Pascal theorem then implies
that the points $A_1$, $A_2$, $B_1$, $B_2$, $C_1$ and $C_2$ lie on
the same conic. \hfill $\square$

\medskip

Let $D_1$, $D_2$, $E_1$, $E_2$, $F_1$ and $F_2$ be the second
intersection points of the conic $\mathcal{C}$ and the lines $A
A_1$, $A A_2$, $B B_1$, $B B_2$, $C C_1$ and $C C_2$,
respectively.

Let $B_3$ be the intersection point of the lines $A_1 C_2$ and
$F_1 D_2$ and $B_4$ be the intersection point of the lines $C_1
A_2$ and $D_1 F_2$. The points $C_3$, $C_4$, $A_3$ and $A_4$ are
defined analogously.

Let $E_3$ be the intersection point of the lines $A_1 C_1$ and
$D_2 F_2$ and $E_4$ be the intersection point of the lines $A_2
C_2$ and $D_1 F_1$. The points $F_3$, $F_4$, $D_3$ and $D_4$ are
defined analogously.

\begin{theorem}\label{t1} The points $B_3$, $B_4$, $E_3$ and $E_4$ lie on the line $C A$,
the points $C_3$, $C_4$, $F_3$ and $F_4$ lie on the line $A B$ and
the points $A_3$, $A_4$, $D_3$ and $D_4$ lie on the line $B C$.
\end{theorem}

\noindent \textbf{Proof:} We shall prove that $B_3$ lie on the
line $C A$.

Let $R$ be the intersection of $A_1 C_2$ and $A C$ and $R'$ the
intersection of the lines $F_1 D_2$ and $A C$.

\begin{figure}[h!h!h!]
\centerline{\includegraphics[width=\textwidth]{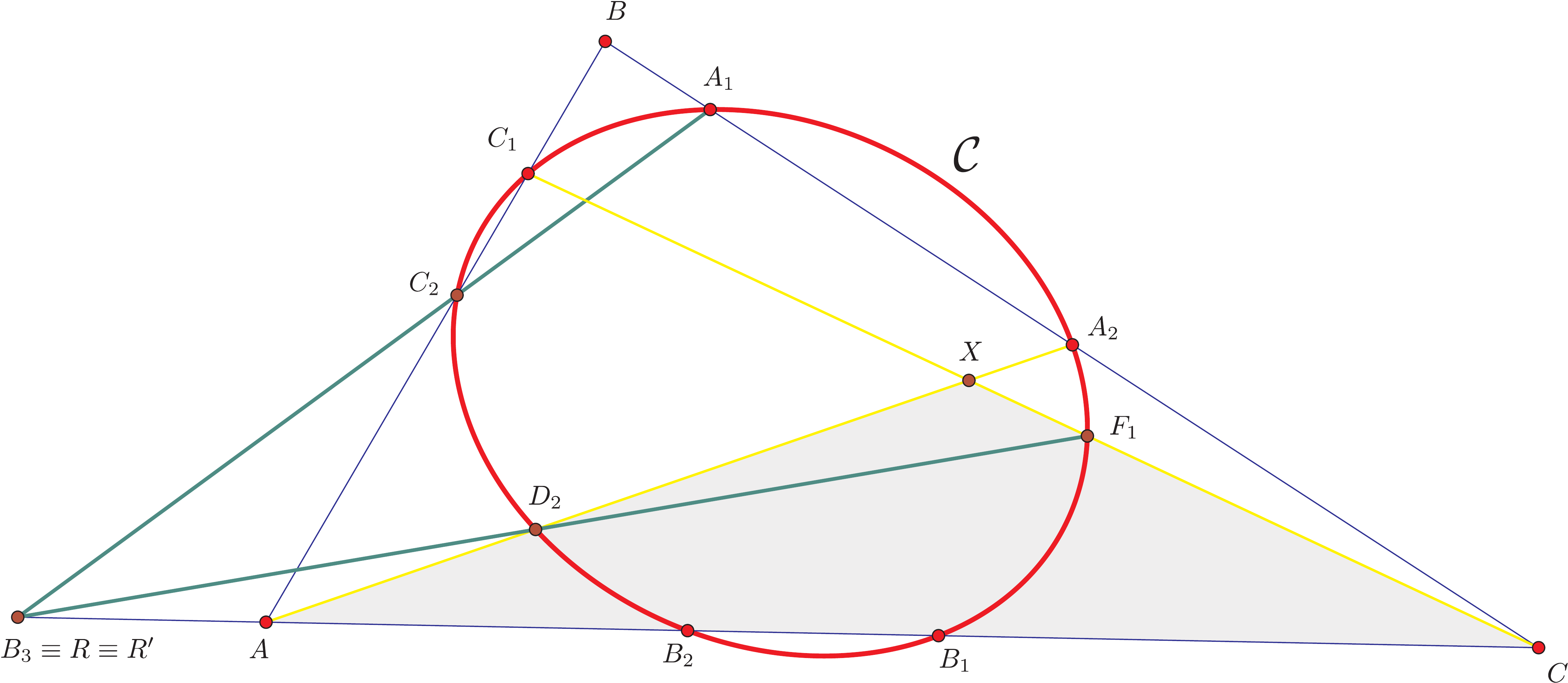}}
\caption{Theorem \ref{t1}} \label{carnot2}
\end{figure}

From the Menelaos theorem for the line $A_1 C_2$, we obtain:
\begin{equation}\label{j1}
\frac{\overrightarrow{C R}}{\overrightarrow{R
A}}=-\frac{\overrightarrow{C_2 B}}{\overrightarrow{A C_2}} \cdot
\frac{ \overrightarrow{A_1 C}}{\overrightarrow{B A_1}}.
\end{equation}

\medskip Let $X$ be the intersection point of the lines $A A_2$
and $C C_1$. The Menelaos theorem for the line $F_1 D_2$ and
$\triangle A X C$ yields:
\begin{equation}\label{j3}
\frac{\overrightarrow{C R'}}{\overrightarrow{R' A}}\cdot
\frac{\overrightarrow{A D_2}}{\overrightarrow{D_2 X}}\cdot
\frac{\overrightarrow{X F_1}}{\overrightarrow{F_1 C}}=-1.
\end{equation}
From the Carnot theorem for the conic $\mathcal{C}$ and $\triangle
A X C$ we obtain:
\begin{equation}\label{j4}
\frac{\overrightarrow{A D_2}}{\overrightarrow{D_2 X}} \cdot
\frac{\overrightarrow{A A_2}}{\overrightarrow{A_2 X}} \cdot
\frac{\overrightarrow{X F_1}}{\overrightarrow{F_1 C}} \cdot
\frac{\overrightarrow{X C_1}}{\overrightarrow{C_1 C}} \cdot
\frac{\overrightarrow{C B_1}}{\overrightarrow{B_1 A}} \cdot
\frac{\overrightarrow{C B_2}}{\overrightarrow{B_2 A}}=1.
\end{equation}
By the Law of Sines we have: $$\overrightarrow{A_2 X} \sin
\sphericalangle A_2 X C= \overrightarrow{A_2 C} \sin
\sphericalangle B C C_1$$ and
$$\overrightarrow{C_1 C} \sin \sphericalangle B C C_1=\overrightarrow{C_1 B} \sin \sphericalangle
\beta.$$ From these two equations one can deduce:
\begin{equation}\label{j5}
\overrightarrow{A_2 X} \cdot \overrightarrow{C_1 C} \sin
\sphericalangle A_2 X C = \overrightarrow{A_2 C} \cdot
\overrightarrow{C_1 B} \sin \sphericalangle \beta.
\end{equation}

Similarly, the following equality holds:
\begin{equation}\label{j6} \overrightarrow{A A_2} \cdot
\overrightarrow{X C_1 } \sin \sphericalangle C_1 X A =
\overrightarrow{B A_2 } \cdot \overrightarrow{A C_1} \sin
\sphericalangle \beta.
\end{equation}

From (\ref{j5}) and (\ref{j6}) (using the equality
$\sphericalangle C_1 X A=\sphericalangle A_2 X C$) we conclude
that:
\begin{equation}\label{j7}
\frac{\overrightarrow{A A_2}}{\overrightarrow{A_2 X}} \cdot
\frac{\overrightarrow{X C_1}}{\overrightarrow{C_1
C}}=\frac{\overrightarrow{B A_2}}{\overrightarrow{A_2
C}}\cdot\frac{\overrightarrow{A C_1}}{\overrightarrow{C_1 B}}.
\end{equation}

Now, from the relations (\ref{j3}), (\ref{j4}) and (\ref{j7}) we
have:
$$\frac{\overrightarrow{C R'}}{\overrightarrow{R' A}}=\frac{\overrightarrow{B A_2}}{\overrightarrow{A_2
C}}\cdot\frac{\overrightarrow{A C_1}}{\overrightarrow{C_1 B}}\cdot
\frac{\overrightarrow{C B_1}}{\overrightarrow{B_1 A}} \cdot
\frac{\overrightarrow{C B_2}}{\overrightarrow{B_2 A}}.$$ But, the
Carnot relation \ref{r} implies $$\frac{\overrightarrow{C
R'}}{\overrightarrow{R' A}}=-\frac{\overrightarrow{C_2
B}}{\overrightarrow{A C_2}} \cdot \frac{ \overrightarrow{A_1
C}}{\overrightarrow{B A_1}},$$ and $R\equiv R'\equiv B_3$.

The proof for the other points is analogous. \hfill $\square$

\begin{figure}[h!h!h!]
\centerline{\includegraphics[width=\textwidth]{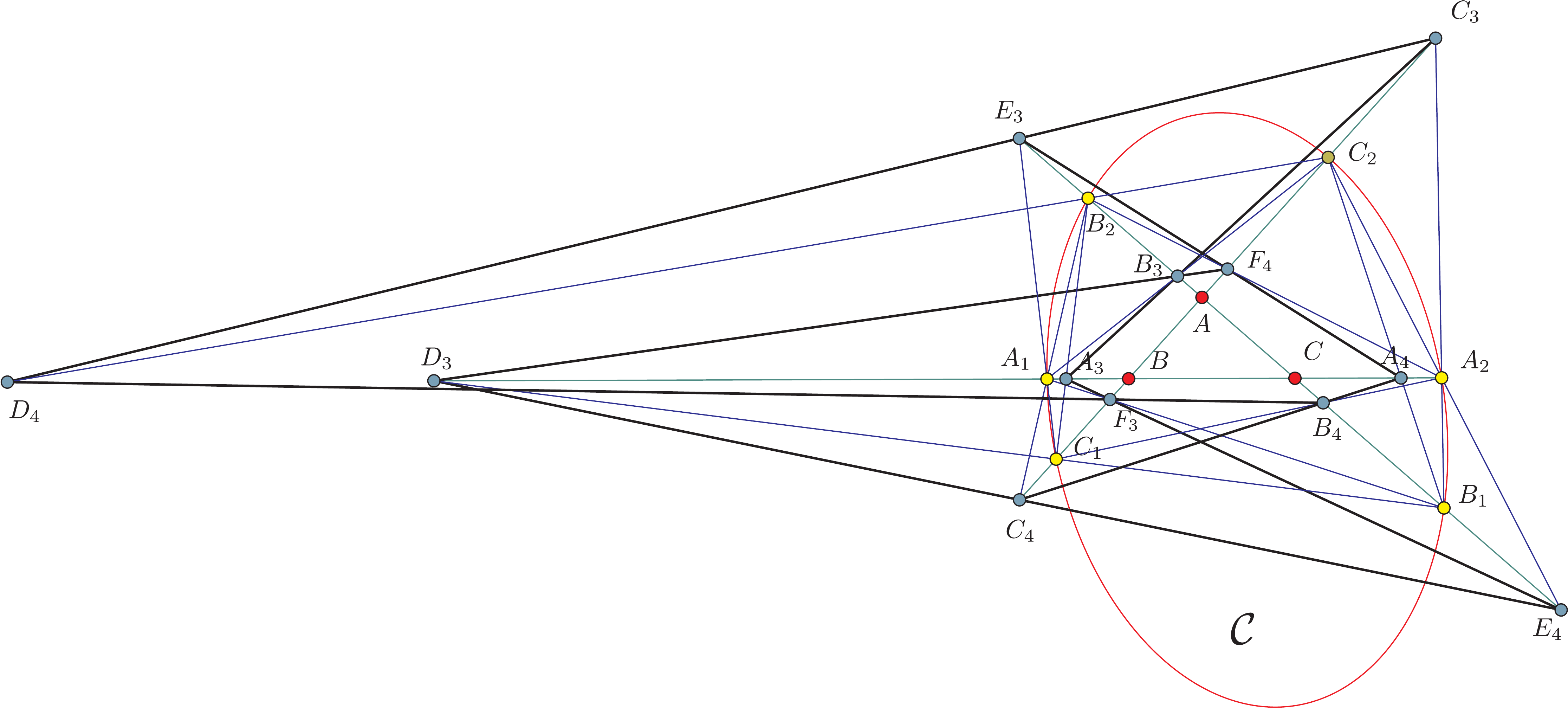}}
\caption{Theorem \ref{t2}} \label{carnot3}
\end{figure}

From Pascal's theorem the following theorem is true (see Figure
\ref{carnot3}):

\begin{theorem} \label{t2} The following 8 triples of points ($A_3$, $B_3$,  $C_3$), ($D_3$, $E_3$,  $C_4$), ($A_3$, $E_4$, $F_3$), ($D_3$, $B_3$,  $F_4$), ($A_4$, $E_3$,  $F_4$),
($D_4$, $E_3$,  $C_3$), ($D_4$, $B_4$,  $F_3$),  and ($A_4$,
$B_4$, $C_4$) are collinear.
\end{theorem}

In the sequel, we encounter the relations of higher order. We use
the theorem of Carnot to prove that certain points in the
configuration lie on the same conic.

\begin{theorem}\label{t4} The points $D_3$, $D_4$, $E_3$, $E_4$, $F_3$ and
$F_4$ lie on the same conic $\mathcal{D}$.
\end{theorem}

\begin{figure}[h!h!h!]
\centerline{\includegraphics[width=0.8 \textwidth]{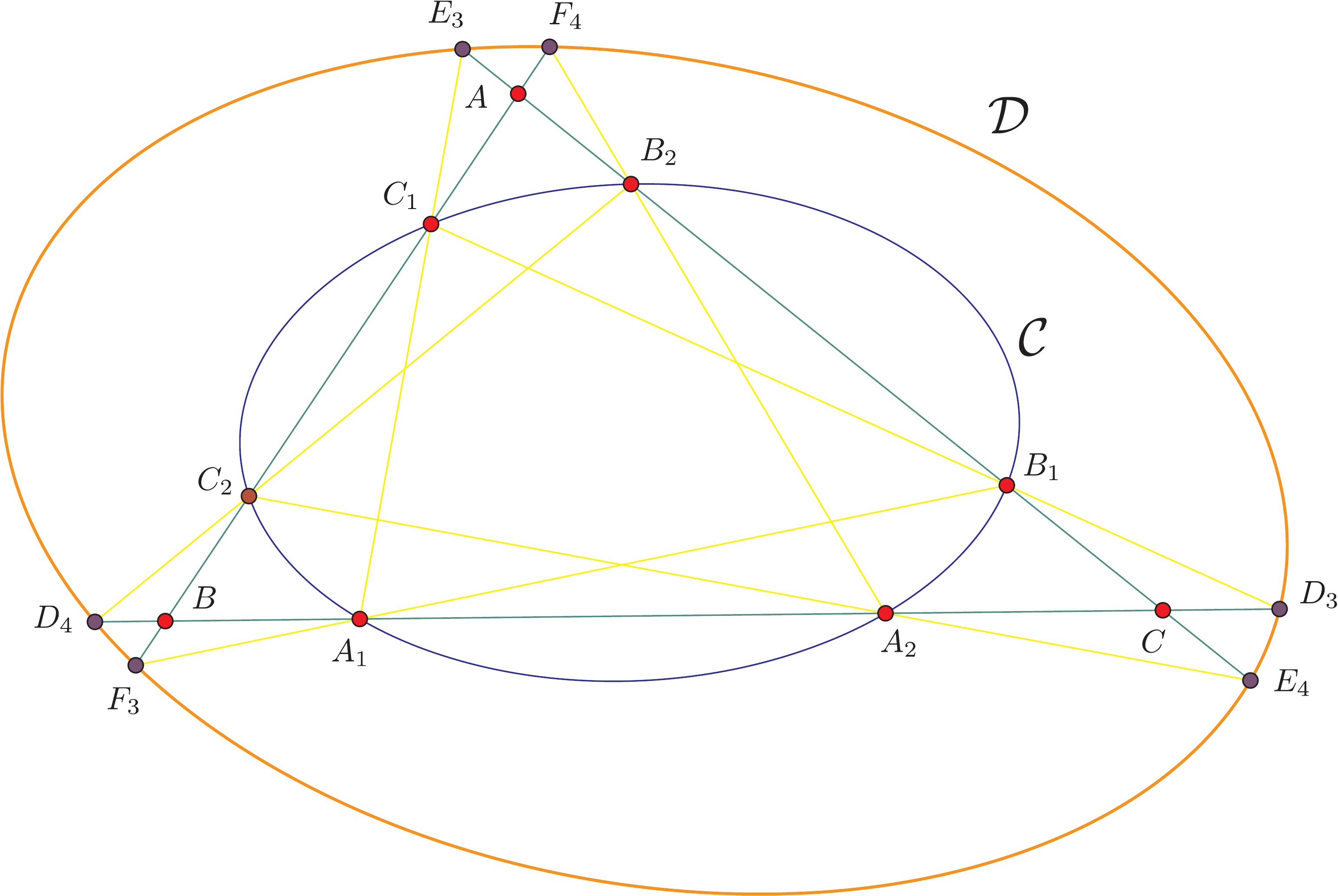}}
\caption{Theorem \ref{t4}} \label{carnot4}
\end{figure}

\noindent \textbf{Proof:} From the proof of Theorem \ref{t1} we
also deduce that: $$\frac{\overrightarrow{C
E_3}}{\overrightarrow{E_3 A}}=-\frac{\overrightarrow{C_1
B}}{\overrightarrow{A C_1}} \cdot \frac{ \overrightarrow{A_1
C}}{\overrightarrow{B A_1}},\, \frac{\overrightarrow{A
F_3}}{\overrightarrow{F_3 B}}=-\frac{\overrightarrow{A_1
C}}{\overrightarrow{B A_1}} \cdot \frac{ \overrightarrow{B_1
C}}{\overrightarrow{C B_1}},\,\frac{\overrightarrow{B
D_3}}{\overrightarrow{D_3 C}}=-\frac{\overrightarrow{B_1
A}}{\overrightarrow{C B_1}} \cdot \frac{ \overrightarrow{C_1
B}}{\overrightarrow{A C_1}},$$
$$\frac{\overrightarrow{C
E_4}}{\overrightarrow{E_4 A}}=-\frac{\overrightarrow{C_2
B}}{\overrightarrow{A C_2}} \cdot \frac{ \overrightarrow{A_1
C}}{\overrightarrow{B A_1}},\,\frac{\overrightarrow{A
F_4}}{\overrightarrow{F_4 B}}=-\frac{\overrightarrow{A_2
C}}{\overrightarrow{B A_2}} \cdot \frac{ \overrightarrow{B_2
C}}{\overrightarrow{C B_2}},\,\frac{\overrightarrow{B
D_4}}{\overrightarrow{D_4 C}}=-\frac{\overrightarrow{B_2
A}}{\overrightarrow{C B_2}} \cdot \frac{ \overrightarrow{C_2
B}}{\overrightarrow{A C_2}}.$$

Then the following holds by \ref{r}:
 $$\frac{\overrightarrow{C
E_3}}{\overrightarrow{E_3 A}}\cdot \frac{\overrightarrow{C
E_4}}{\overrightarrow{E_4 A}}\cdot \frac{\overrightarrow{A
F_3}}{\overrightarrow{F_3 B}}\cdot \frac{\overrightarrow{A
F_4}}{\overrightarrow{F_4 B}}\cdot \frac{\overrightarrow{A
F_3}}{\overrightarrow{F_3 B}}\cdot \frac{\overrightarrow{B
D_4}}{\overrightarrow{D_4 C}}=1.$$ By the converse of Carnot's
theorem, the points $D_3$, $D_4$, $E_3$, $E_4$, $F_3$ and $F_4$
lie on the same conic. \hfill $\square$

In the same fashion we prove that:
\begin{theorem}\label{t5} The following 4 sextuples of the points $(D_3, D_4, E_3, E, 4, F_3, F_4)$,\\ $(A_3, A_4, B_3, B_4, F_3, F_4)$, $(A_3, A_4, E_3, E_4, C_3, C_4)$ and $(D_3, D_4, B_3, B_4, C_3, C_4)$ are the
sextuples of the points lying on the same conic.
\end{theorem}

\begin{figure}[h!h!h!]
\centerline{\includegraphics[width=\textwidth]{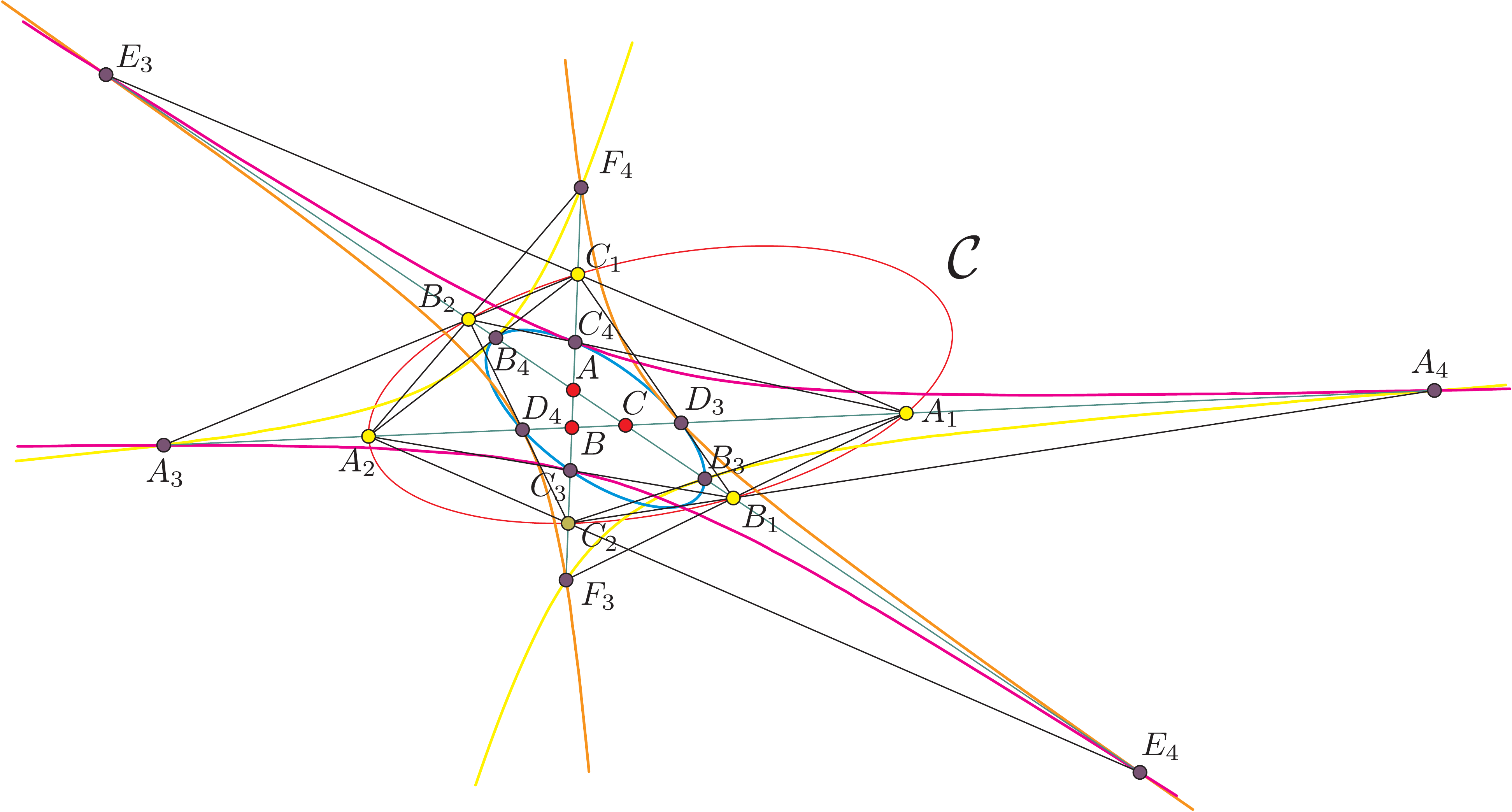}}
\caption{Theorem \ref{t5}} \label{carnot5}
\end{figure}

\section{Bradley's Theorem}

In this section we give an elementary proof of the Bradley's
conjecture \cite{Bradley}. The first proof, given by Zolt\'{a}n
Szilasi in \cite{Szilasi}, used barycentric coordinates. We use
different approach and prove several other interesting things
about Carnot's configuration.

Let $X_1$ be the intersection points of the lines $A A_1$ and $B
B_1$, $X_2$ of $B B_1$ and $C C_1$ and $X_3$ of $C C_1$ and $A
A_1$. Let $Y_1$ be the intersection points of the lines $A A_2$
and $B B_2$, $Y_2$ of $B B_2$ and $C C_2$ and $Y_3$ of $C C_2$ and
$A A_2$.

Define $T_2$ as the intersection point of the lines $X_1 Y_3$ and
$X_3 Y_1$. The points $T_3$ and $T_1$ are defined analogously.

\begin{theorem}
   $T_2$ lies on the line $B C$, $T_3$ on $C A$ and $T_1$
    on $A B$.
\end{theorem}

\noindent \textbf{Proof:} Let $T'$ be the intersection point of
the lines $X_3 Y_1$ and $B C$ and let $T''$ be the intersection
point of the lines $X_1 Y_3$ and $B C$.

By the Menelaos theorem applied at $\triangle A B A_1$ and the
line $B C_2$ we obtain:
 $$\frac{\overrightarrow{A
X_3}}{\overrightarrow{X_3 A_1}}=-\frac{\overrightarrow{C
B}}{\overrightarrow{A_1 C}} \cdot \frac{ \overrightarrow{A
C_1}}{\overrightarrow{C_1 B}}.$$ The same reasoning for $\triangle
A C A_2$ and the line $C C_1$ we obtain:
$$\frac{\overrightarrow{A_2 Y_1}}{\overrightarrow{Y_1
A}}=-\frac{\overrightarrow{B_2 C}}{\overrightarrow{A B_2}} \cdot
\frac{ \overrightarrow{B A_2}}{\overrightarrow{C B}}.$$ Then from
the Menelaos theorem for $\triangle A A_1 A_2$ and the line $X_1
Y_3$  we get:  \begin{equation}\label{t'}\frac{\overrightarrow{A_1
T'}}{\overrightarrow{T' A_2}}=-\frac{\overrightarrow{A
B_2}}{\overrightarrow{A C_1}} \cdot \frac{ \overrightarrow{A_1
C}}{\overrightarrow{B_2 C}} \cdot \frac{ \overrightarrow{C_1
B}}{\overrightarrow{B A_2}}.\end{equation} In the same fashion we
prove that: \begin{equation}\label{t"}\frac{\overrightarrow{A_1
T''}}{\overrightarrow{T'' A_2}}=-\frac{\overrightarrow{A
C_2}}{\overrightarrow{A B_1}} \cdot \frac{ \overrightarrow{B_1
C}}{\overrightarrow{A_2 C}} \cdot \frac{ \overrightarrow{B
A_1}}{\overrightarrow{C_2 B}}.\end{equation}

By the relation (\ref{r}) we conclude that:
$$\frac{\overrightarrow{A_1 T'}}{\overrightarrow{T' A_2}}=\frac{\overrightarrow{A_1
T''}}{\overrightarrow{T'' A_2}}, $$ so $T'\equiv T''\equiv T_2$.

For the points $T_1$ and $T_3$ the proof is analogous. \hfill
$\square$

Since the points $T_2$, $B$ and $C$ are collinear, by the converse
of Pascal's theorem for the hexagon $X_3 Y_1 Y_2 Y_3 X_1 X_2$ we
get (see Figure \ref{bradley1s}):

\begin{figure}[h!h!h!]
\centerline{\includegraphics[width=\textwidth]{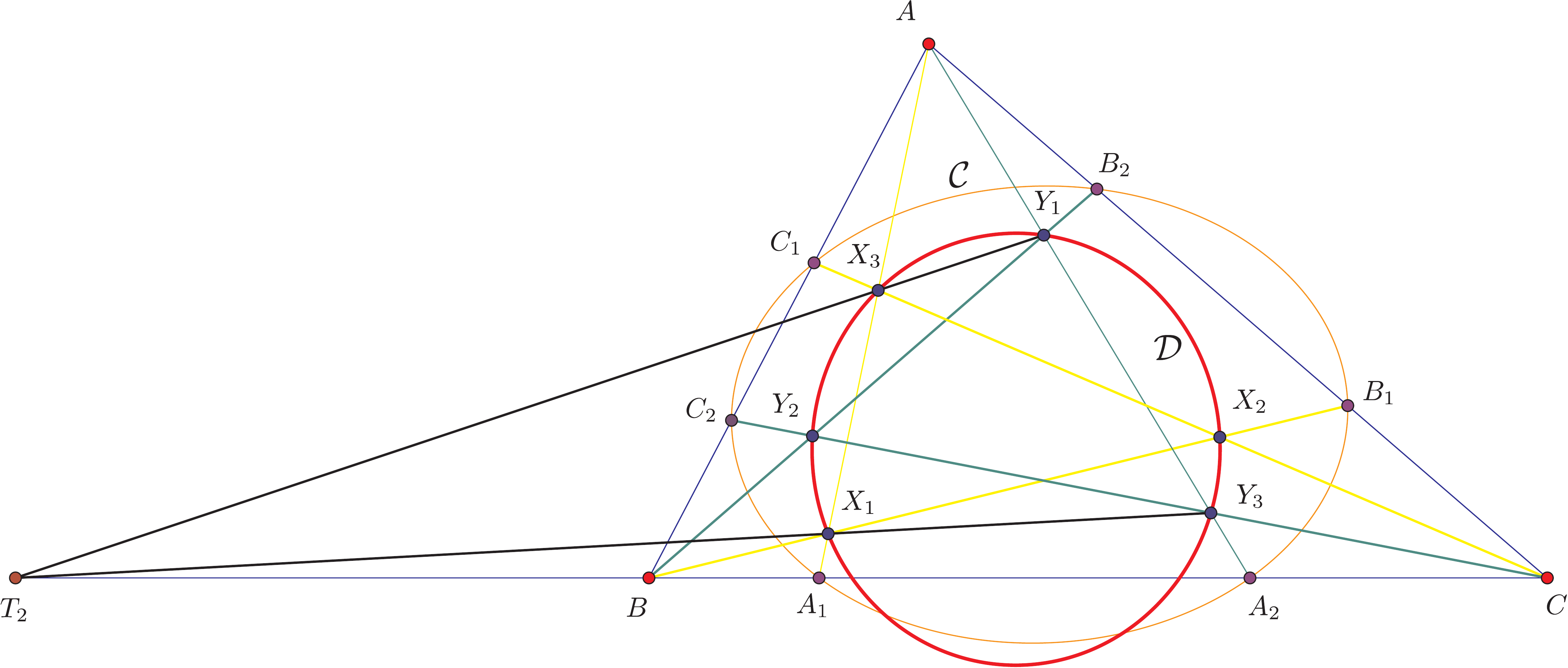}}
\caption{Corollary \ref{bradley1}} \label{bradley1s}
\end{figure}

\begin{corollary}\label{bradley1} The points $X_1$, $X_2$, $X_3$, $Y_1$, $Y_2$ and
$Y_3$ lie on the same conic.
\end{corollary}

\begin{figure}[h!h!h!]
\centerline{\includegraphics[width=0.7\textwidth]{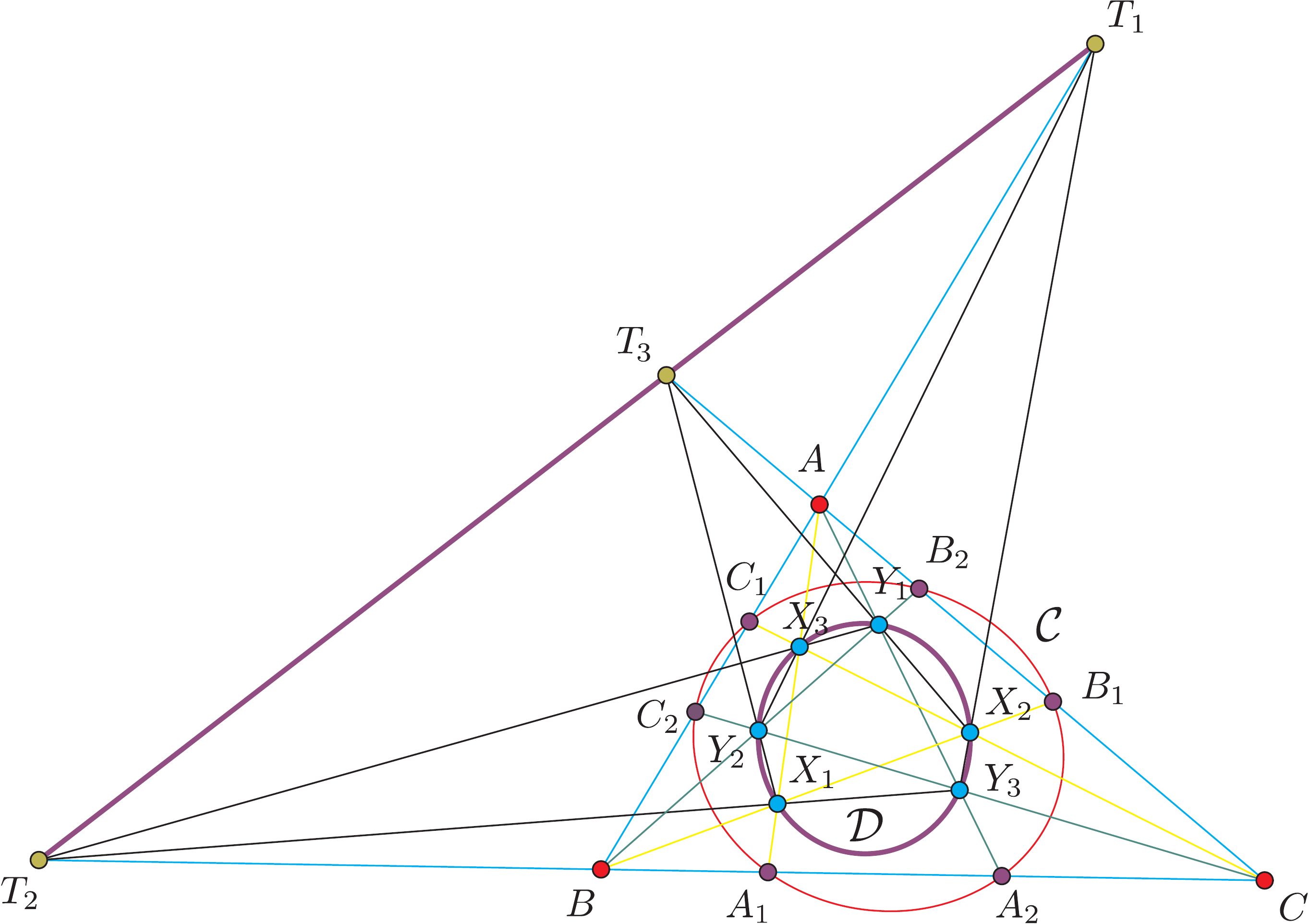}}
\caption{Corollary \ref{bradley2}} \label{bradley2s}
\end{figure}

An immediate consequence of this fact is (see Figure
\ref{bradley2s}):

\begin{corollary} \label{bradley2} The points $T_1$, $T_2$ and $T_3$ lie on the same line.
\end{corollary}

\begin{figure}[h!h!h!]
\centerline{\includegraphics[width=0.75\textwidth]{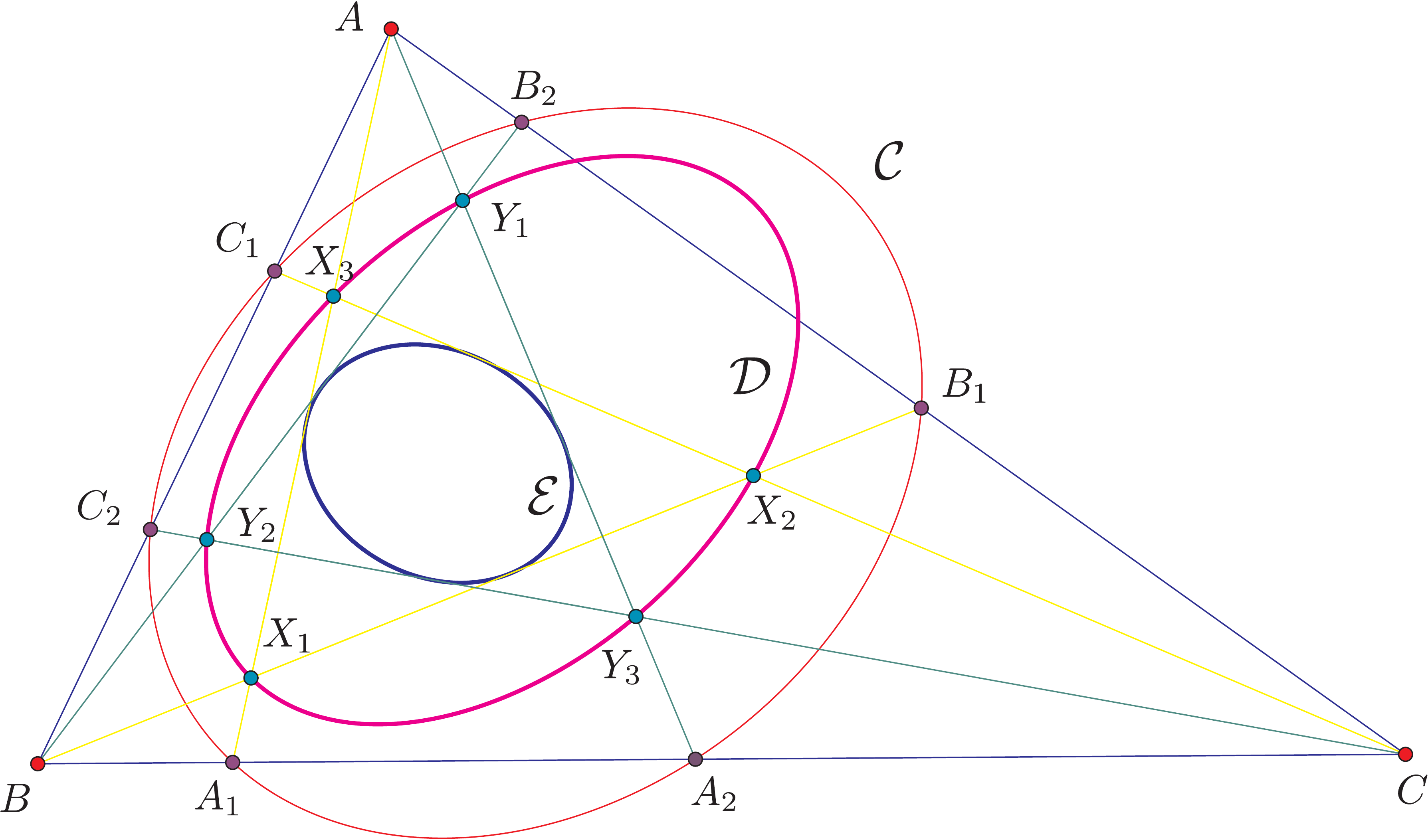}}
\caption{Bradley's theorem} \label{bradley3s}
\end{figure}

Bradley's theorem \ref{glavna},  directly follows from Corollary
\ref{bradley1} and the Poncelet triangle theorem \cite[Theorem 5,
p.184-185]{Pra}, see Figure \ref{bradley3s}.

\section{Proof of Theorem \ref{glavnag}}

In this section we give the proof of Theorem \ref{glavnag}. The
proof illustrates a nice application of the Menelaus and the
Carnot theorems.

\noindent \textbf{Proof:} We prove that the points $23$, $24$,
$32$, $34$, $42$, $43$ lie on conic $\mathcal{C}_1$. The proof for
other points is analogous.

\begin{figure}[h!h!h!]
\centerline{\includegraphics[width=0.8\textwidth]{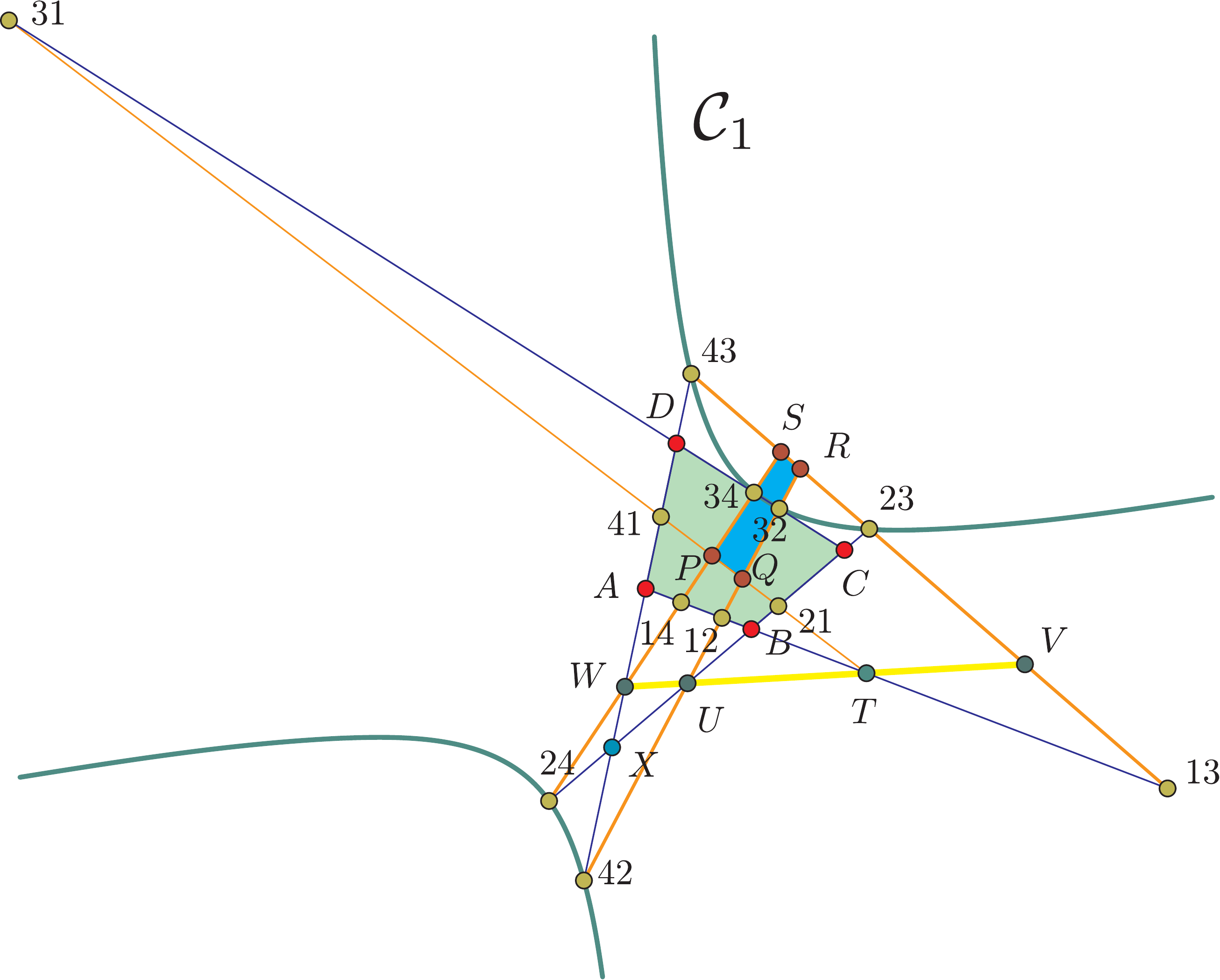}}
\caption{Bradley's theorem about quadrilaterals} \label{bradleyqs}
\end{figure}

Let $X$ be the intersection point of the lines $A D$ and $B C$. We
apply the Menelaus theorem for $\triangle X D C$ and the lines $S
W$ , $R U$, $S V$ and $V W$ and get:

\begin{equation}\label{g1}
\frac{\overrightarrow{X W}}{\overrightarrow{W D}}\cdot
\frac{\overrightarrow{D (34)}}{\overrightarrow{(34) C}}\cdot
\frac{\overrightarrow{C (24)}}{\overrightarrow{(24) X}}=-1,
\end{equation}
\begin{equation}\label{g2}
\frac{\overrightarrow{X (41)}}{\overrightarrow{(41) D}}\cdot
\frac{\overrightarrow{D (32)}}{\overrightarrow{(32) C}}\cdot
\frac{\overrightarrow{C U}}{\overrightarrow{(U X}}=-1,
\end{equation}
\begin{equation}\label{g3}
\frac{\overrightarrow{X (43)}}{\overrightarrow{(43)
D}}\cdot\frac{\overrightarrow{D V}}{\overrightarrow{V C}}\cdot
\frac{\overrightarrow{C (23)}}{\overrightarrow{(23) X}}=-1,
\end{equation}
\begin{equation}\label{g4}
\frac{\overrightarrow{D W}}{\overrightarrow{W X}}\cdot
\frac{\overrightarrow{X U}}{\overrightarrow{U C}}\cdot
\frac{\overrightarrow{C V}}{\overrightarrow{V D}}=-1.
\end{equation}

After multiplication of (\ref{g1}), (\ref{g2}), (\ref{g3}) and
(\ref{g4}), we obtain: $$ \frac{\overrightarrow{D
(34)}}{\overrightarrow{(34) C}}\cdot \frac{\overrightarrow{C
(24)}}{\overrightarrow{(24) X}}\cdot \frac{\overrightarrow{X
(41)}}{\overrightarrow{(41) D}}\cdot \frac{\overrightarrow{D
(32)}}{\overrightarrow{(32) C}}\cdot \frac{\overrightarrow{X
(43)}}{\overrightarrow{(43) D}}\cdot \frac{\overrightarrow{C
(23)}}{\overrightarrow{(23) X}}=1.$$ From the converse of Carnot's
theorem it follows that the points $23$, $24$, $32$, $34$, $42$,
$43$ lie on the same conic. \hfill $\square$

\begin{center}\textmd{Acknowledgements }
\end{center}

\medskip This research is done during my stay in Switzerland. The
author wishes to thank his friends the Hajdin family: Katarina,
Rade, Nikola, Luka and Matija for generous hospitality and
support.

\medskip

{\small \DJ{}OR\DJ{}E BARALI\'{C}, Mathematical Institute SASA,
Kneza Mihaila 36, p.p.\ 367, 11001 Belgrade, Serbia

E-mail address: djbaralic@mi.sanu.ac.rs


\begin{thebibliography}{99} 


\bibitem{Bradley} C. Bradley, Problems requiring proofs,\\
avaiable http://people.bath.ac.uk/masgcs/Article182.pdf, 2011


\bibitem{Hatt} J. L. S. Hatton, \textit{The Principles of Projective Geometry
Applied to the Straight Line and Conic}, Cambridge University
Press 1913.


\bibitem{Oste} A. Ostermann and G. Wanner, \textit{Geometry by Its History}, Springer,
2012.


\bibitem{Pra}  V. V. Prasolov and V. M. Tikhomirov,
\textit{Geometry}, American Mathematical Society, 2001.


\bibitem{Gerb}  J. Richter-Gebert,
\textit{Perspectives on Projective Geometry}, Springer-Verlag
Berlin Heidelberg, 2011.

\bibitem{JRG} J. Richter-Gebert, \textit{ Meditations on Ceva's Theorem}, In The
Coxeter Legacy: Reflections and Projections (Eds. Chandler Davis
\& Eric Ellers, American Mathematical Society, Fields Institute),
227--254, 2006.

\bibitem{Szilasi} Z. Szilasi, Two applications of the theorem of
Carnot, Annales Mathematicae et Informaticae \textbf{40} (2012),
135-144.


\end{thebibliography}
\end{document}